\documentclass[conference]{IEEEtran}
\usepackage{amsmath,amssymb,amsthm,multicol}
\usepackage{url}

\begin{document}

\newtheorem{thm}{Theorem}[section]
\newtheorem{cor}{Corollary}[section]
\newtheorem{lemm}{Lemma}[section]
\newtheorem{conj}{Conjecture}[section]
\theoremstyle{definition}
\newtheorem{defn}{Definition}[section]
\newtheorem{exmpl}{Example}[section]
\newtheorem{remrk}{Remark}[section]

\title{Solving the Odd Perfect Number Problem: \\Some New Approaches}
\author{\IEEEauthorblockN{Jose Arnaldo B. Dris}
\IEEEauthorblockA{Master of Science in Mathematics Graduate (2008) \\
De La Salle University, 2401 Taft Avenue, 1004 Manila, Philippines \\
jabdris@yahoo.com.ph \\}
}

\maketitle

\begin{abstract}
\boldmath
A conjecture predicting an injective and surjective mapping $X = \displaystyle\frac{\sigma(p^k)}{p^k}, Y = \displaystyle\frac{\sigma(m^2)}{m^2}$ between OPNs $N = {p^k}{m^2}$ (with Euler factor $p^k$) and rational points on the hyperbolic arc $XY = 2$ with $1 < X < 1.25 < 1.6 < Y < 2$ and $2.85 < X + Y < 3$, is disproved.  We will show that if an OPN $N$ has the form above, then $p^k < {\frac{2}{3}}{m^2}$.  We then give a somewhat weaker corollary to this last result ($m^2 - p^k \ge 8$) and give possible improvements along these lines.  We will also attempt to prove a conjectured improvement to $p^k < m$ by observing that $\displaystyle\frac{\sigma(p^k)}{m} \ne 1$ and $\displaystyle\frac{\sigma(p^k)}{m} \ne \displaystyle\frac{\sigma(m)}{p^k}$ in all cases.  Lastly, we also prove the following generalization:  If $N = \displaystyle\prod_{i = 1}^{r}{{p_i}^{\alpha_i}}$ is the canonical factorization of an OPN $N$, then $\displaystyle\sigma({p_i}^{\alpha_i}) \le \displaystyle\frac{2}{3}\displaystyle\frac{N}{{p_i}^{\alpha_i}}$ for all $ i$.  This gives rise to the inequality $N^{2 - r} \le (\frac{1}{3}){\left(\frac{2}{3}\right)}^{r - 1}$, which is true for all $r$, where $r = \omega(N)$ is the number of distinct prime factors of $N$.
\end{abstract}


\section{INTRODUCTION}
A perfect number is a positive integer $N$ such that the sum of all the positive divisors of $N$ equals $2N$, denoted by $\sigma(N) = 2N$. The question of the existence of odd perfect numbers (OPNs) is one of the longest unsolved problems of number theory.  This paper gives some new and nontraditional attempts and approaches to solving the Odd Perfect Number (OPN) Problem. \\

Hereinafter, we shall let $N = {p^k}{m^2}$ denote an OPN with Euler (prime-power) factor $p^k$ (with $p \equiv k \equiv 1 \pmod 4$ and $\gcd(p, m) = 1$), assuming at least one such number exists.

\section{SOME PREPARATORY LEMMAS}
The following lemmas will be very useful later on:

\begin{lemm}\label{Lemma1}
\hspace{0.01in} \\
$0 < \rho_2 \le \frac{2}{3} < 1 < \rho_1 < \frac{5}{4} < \frac{8}{5} < \mu_1 < 2 < 3 \le \mu_2$ where:
$$\rho_1 = \displaystyle\frac{\sigma(p^k)}{p^k}$$
$$\rho_2 = \displaystyle\frac{\sigma(p^k)}{m^2}$$
\hspace{0.01in} \\
$$\mu_1  = \displaystyle\frac{\sigma(m^2)}{m^2}$$
$$\mu_2 = \displaystyle\frac{\sigma(m^2)}{p^k}$$
\end{lemm}

\begin{IEEEproof}
We refer the interested reader to the author's master's thesis completed in August of 2008 \cite{Dris}.
\end{IEEEproof}

\begin{lemm}\label{Lemma2}
$\frac{57}{20} < \rho_1 + \mu_1 < 3 < \frac{11}{3} \le \rho_2 + \mu_2$
\end{lemm}

\begin{IEEEproof}
Again, the interested reader is referred to \cite{Dris}, where it is shown that "any further improvement to the lower bound of $\frac{57}{20}$ for $\rho_1 + \mu_1$ would be equivalent to showing that there are no OPNs of the form $5{m^2}$ which would be a very major result. Likewise, any further improvement on the upper bound of 3 would have similar implications for all arbitrarily large primes and thus would be a very major result."  (These assertions, which are originally Joshua Zelinsky's, were readily verified by the author using Mathematica.) 
\end{IEEEproof}

\begin{lemm}\label{Lemma3}
Define $\rho_3 = \displaystyle\frac{\sigma(p^k)}{m}$ and $\mu_3 = \displaystyle\frac{\sigma(m)}{p^k}$.  Then $\rho_3 \ne \mu_3$, and the following statements hold:
\begin{itemize}
\item{If $\rho_3 < 1$, then $p^k < m$.}
\item{Suppose that $1 < \rho_3$.
\begin{itemize}
\item{If $\rho_3 < \mu_3$, then $\frac{4}{5}{m} < p^k < \sqrt{2}{m}$.}
\item{If $\mu_3 < \rho_3$, then $m < p^k$.}
\end{itemize}
}
\end{itemize}
\end{lemm}

\begin{IEEEproof}
The interested reader is again referred to \cite{Dris}.  The crucial part of the argument is in showing that $\displaystyle\frac{\sigma(p^k)}{p^k} < \displaystyle\frac{\sigma(m)}{m}$.
\end{IEEEproof}

\section{MAIN RESULTS}
First, we prove that a conjectured one-to-one correspondence is actually both not surjective and not injective.

\begin{conj}\label{Conjecture1}
For each $N = {p^k}{m^2}$ an OPN with \\$N > {10}^{300}$, there corresponds exactly one ordered pair of rational numbers $\left(\displaystyle\frac{\sigma(p^k)}{p^k}, \displaystyle\frac{\sigma(m^2)}{m^2}\right)$ lying in the region \\$1 < \displaystyle\frac{\sigma(p^k)}{p^k} < \displaystyle\frac{5}{4}$, $\displaystyle\frac{8}{5} < \displaystyle\frac{\sigma(m^2)}{m^2} < 2$, and \\ $\displaystyle\frac{57}{20} < \displaystyle\frac{\sigma(p^k)}{p^k} + \displaystyle\frac{\sigma(m^2)}{m^2} < 3$, and vice-versa.
\end{conj}

\begin{IEEEproof}
(Note that this is actually a refutation of the conjecture.)  First, we note that the equation $\displaystyle\frac{\sigma(x)}{x} = \displaystyle\frac{\sigma(p^k)}{p^k}$ has the sole solution $x = p^k$ since prime powers are solitary.  Next, let $X = \displaystyle\frac{\sigma(p^k)}{p^k}$ and $Y = \displaystyle\frac{\sigma(m^2)}{m^2}$.  It is straightforward to observe that, since the abundancy index is an arithmetic function, then for each $N = {p^k}{m^2}$ an OPN (with $N > {10}^{300}$), there corresponds exactly one ordered pair of rational numbers $(X, Y)$ lying in the hyperbolic arc $XY = 2$ bounded as follows: $1 < X < 1.25$, $1.6 < Y < 2$, and $2.85 < X + Y < 3$.  (Note that these bounds are the same ones obtained in Lemma \ref{Lemma1} and Lemma \ref{Lemma2}.) \\

We now disprove the backward direction of the conjecture.  We do this by showing that the mapping $X = \displaystyle\frac{\sigma(p^k)}{p^k}$ and $Y = \displaystyle\frac{\sigma(m^2)}{m^2}$ is neither surjective nor injective in the specified region.

\emph{$(X, Y)$ is not surjective.} We prove this claim by producing a rational point $(X_0, Y_0)$ lying in the specified region, and which satisfies $X_0 = \displaystyle\frac{\sigma(pq)}{pq}$ where $p$ and $q$ are primes satisfying $5 < p < q$.  Notice that
$$1 < X_0 = \displaystyle\frac{(p + 1)(q + 1)}{pq} = \left(1 + \displaystyle\frac{1}{p}\right)\left(1 + \displaystyle\frac{1}{q}\right) \le \displaystyle\frac{8}{7}\displaystyle\frac{12}{11} = \displaystyle\frac{96}{77}$$
where $\displaystyle\frac{96}{77} < 1.2468 < 1.25$.  Now, by setting $Y_0 = \displaystyle\frac{2}{X_0}$, the other two inequalities for $Y_0$ and $X_0 + Y_0$ would follow.  Thus, we now have a rational point $(X_0, Y_0)$ in the specified region, and which, by a brief inspection, satisfies $X_0 \ne \displaystyle\frac{\sigma(p^k)}{p^k}$ for all primes $p$ and positive integers $k$ (since prime powers are solitary). Consequently, the mapping defined in the backward direction of the conjecture is not surjective. \\

\emph{Remark.  Since the mapping is not onto, there are rational points in the specified region which do not correspond to any OPN.} \\

\emph{$(X, Y)$ is not injective.} It suffices to construct two distinct OPNs $N_1 = {{p_1}^{k_1}}{m_1}^2$ and $N_2 = {{p_2}^{k_2}}{m_2}^2$ that correspond to the same rational point $(X, Y)$.  Since it cannot be the case that ${p_1}^{k_1} \ne {p_2}^{k_2}$, ${m_1}^2 = {m_2}^2$, we consider the scenario ${p_1}^{k_1} = {p_2}^{k_2}$, ${m_1}^2 \ne {m_2}^2$.  Thus, we want to produce a pair $(m_1, m_2)$ satisfying $\displaystyle\frac{\sigma({m_1}^2)}{{m_1}^2} = \displaystyle\frac{\sigma({m_2}^2)}{{m_2}^2}$. (A computer check by a foreign professor using Maple produced no examples for this equation in the range $1 \le m_1 < m_2 \le 300000$.  But then again, in pure mathematics, absence of evidence is not evidence of absence.)  Now, from the inequalities ${p^k} < {m^2}$ and $N = {p^k}{m^2} > {10}^{300}$, we get $m^2 > {10}^{150}$.  A nonconstructive approach to finding a solution to $\displaystyle\frac{\sigma({m_1}^2)}{{m_1}^2} = \displaystyle\frac{\sigma({m_2}^2)}{{m_2}^2}$ would then be to consider ${10}^{150} < {m_1}^2 < {m_2}^2$ and Erdos' result that ``\emph{The number of solutions of $\displaystyle\frac{\sigma(a)}{a} = \displaystyle\frac{\sigma(b)}{b}$ satisfying $a < b \le x$ equals $Cx + o(x)$ where $C \ge \displaystyle\frac{8}{147}$.}'' (\cite{Erdos}, \cite{Anderson})  (Note that $C$ here is the same as the (natural) density of friendly integers.)  Given Erdos' result then, this means that eventually, as $m_2 \rightarrow \infty$, there will be at least ${{10}^{150}}\displaystyle\frac{8}{147}$ solutions $(m_1, m_2)$ to $\displaystyle\frac{\sigma({m_1}^2)}{{m_1}^2} = \displaystyle\frac{\sigma({m_2}^2)}{{m_2}^2}$, a number which is obviously greater than $1$.  This finding, though nonconstructive, still proves that the mapping defined in the backward direction of the conjecture is not injective.   
\end{IEEEproof}

Next, we show how the components $p^k$ and $m^2$ of an OPN are related.

\begin{thm}\label{Theorem1}
${p^k} < \displaystyle\frac{2}{3}{m^2}$
\end{thm}

\begin{IEEEproof}
This theorem follows from the inequalities $\rho_1 > 1$ and $\rho_2 \le \displaystyle\frac{2}{3}$.
\end{IEEEproof}

A somewhat weaker result than Theorem \ref{Theorem1} is the following:

\begin{thm}\label{Theorem2}
${m^2} - {p^k} \ge 8$ 
\end{thm}

\begin{IEEEproof}
From Theorem \ref{Theorem1}, we know that ${m^2} - {p^k} > \displaystyle\frac{p^k}{2}$.  So in particular, we are sure that ${m^2} - {p^k} > 0$.  But $m$ odd implies that $m^2 \equiv 1 \pmod 4$, and we also know that $p^k \equiv 1 \pmod 4$.  Thus, ${m^2} - {p^k} \equiv 0 \pmod 4$, which is equivalent to saying that $4|({m^2} - {p^k})$.  Since ${m^2} - {p^k} > 0$, then this implies that ${m^2} - {p^k} \ge 4$.  Suppose ${m^2} - {p^k} = 4$. \\

Then ${p^k} = {m^2} - 4 = (m + 2)(m - 2)$.  Hence, we have the simultaneous equations ${p^{k - x}} = m + 2$ and ${p^x} = m - 2$ where $k \ge 2x + 1$.  Consequently, we have ${p^{k - x}} - {p^x} = 4$, which implies that ${p^x}({p^{k - 2x}} - 1) = 4$ where $k - 2x$ is odd.  Since $(p - 1)|({p^y} - 1)$ $\forall y \ge 1$, this last equation implies that $(p - 1) | 4$.  Likewise, the congruence $p \equiv 1 \pmod 4$ implies that $4 | (p - 1)$.  These two divisibility relations imply that $p - 1 = 4$, or $p = 5$.  Hence, ${5^x}(5^{k - 2x} - 1) = 4$.  Since $5$ does not divide $4$, $x = 0$ and thus ${5^k} - 1 = 4$, which means that $k = 1$. \\

Therefore, $p = 5, k = 1, x = 0$.  Consequently, $p^{k - x} = 5^{1 - 0} = 5 = m + 2$ and ${p^x} = 5^0 = 1 = m - 2$.  Either way, we have $m = 3$. \\

All of these computations imply that $N = {p^k}{m^2} = {5^1}{3^2} = 45$ is an OPN.  But this contradicts the fact that $\displaystyle\frac{\sigma(N)}{N} = \displaystyle\frac{26}{15} < 2$ (i.e., $N = 45$ is deficient). \\

Thus, ${m^2} - {p^k} \ge 8$.
\end{IEEEproof}

Lastly, we prove the following generalization to the inequality $\sigma(p^k) \le {\displaystyle\frac{2}{3}}{m^2} = {\displaystyle\frac{2}{3}}{\displaystyle\frac{N}{p^k}}$. \\

\begin{thm}\label{Theorem3}
Let $N = \displaystyle\prod_{i = 1}^{\omega(N)}{{p_i}^{\alpha_i}}$ be the canonical factorization of an OPN $N$, where $p_1 < p_2 < \ldots < p_t$ are primes, $t = \omega(N)$ and $\alpha_i > 0$ for all $i$.  Then
$\sigma({p_i}^{\alpha_i}) \le \displaystyle\frac{2}{3}\displaystyle\frac{N}{{p_i}^{\alpha_i}}$
for all $i$.
\end{thm}

\begin{IEEEproof}
Let an OPN be given in the form $N = {p_i}^{\alpha_i}M$ for a particular $i$.  Since ${p_i}^{\alpha_i}||N$ and $N$ is an OPN, then $\sigma({p_i}^{\alpha_i})\sigma(M) = 2{p_i}^{\alpha_i}M$.  Since ${p_i}^{\alpha_i}$ and $\sigma({p_i}^{\alpha_i})$ are always relatively prime, we know that ${p_i}^{\alpha_i} | \sigma(M)$ and we have $\sigma(M) = h{{p_i}^{\alpha_i}}$ for some positive integer $h$.  Assume $h$ = 1.  Then $\sigma(M) = {p_i}^{\alpha_i}$, forcing $\sigma({p_i}^{\alpha_i}) = 2M$.  Since $N$ is an OPN, ${p_i}^{\alpha_i}$ is odd, whereupon we have an odd $\alpha_i$ by considering parity conditions from the last equation.  But this means that ${p_i}^{\alpha_i}$ is the Euler's factor of $N$, and we have ${p_i}^{\alpha_i} = p^k$ and $M = m^2$.  Consequently, $\sigma(m^2) = \sigma(M) = {p_i}^{\alpha_i} = p^k$, which contradicts the fact that $\mu_2 \ge 3$.  Now suppose that $h = 2$.  Then we have the equations $\sigma(M) = 2{{p_i}^{\alpha_i}}$ and $\sigma({p_i}^{\alpha_i}) = M$.  (Note that, since $M$ is odd, $\alpha_i$ must be even.)  Applying the $\sigma$ function to both sides of the last equation, we get $\sigma(\sigma({p_i}^{\alpha_i})) = \sigma(M) = 2{{p_i}^{\alpha_i}}$, which means that ${p_i}^{\alpha_i}$ is an odd superperfect number.  But Kanold \cite{Suryanarayana1} showed that odd superperfect numbers must be perfect squares (no contradiction at this point, since $\alpha_i$ is even), and Suryanarayana \cite{Suryanarayana2} showed in 1973 that ``\emph{There is no odd perfect number of the form $p^{2\alpha}$}'' (where $p$ is prime).  Thus, $h = \displaystyle\frac{\sigma(M)}{{p_i}^{\alpha_i}} \ge 3$, whereupon we have the result $\sigma({p_i}^{\alpha_i}) \le \displaystyle\frac{2}{3}M = \displaystyle\frac{2}{3}\displaystyle\frac{N}{{p_i}^{\alpha_i}}$ for the chosen $i$.  Since $i$ was arbitrary, we have proved our claim in this theorem.
\end{IEEEproof}

The following corollary is a direct consequence of Theorem \ref{Theorem3}.

\begin{cor}\label{Corollary1}
Let $N$ be an OPN with $r = \omega(N)$ distinct prime factors.  Then $N^{2 - r} \le {\left(\displaystyle\frac{1}{3}\right)}{\left(\displaystyle\frac{2}{3}\right)}^{r - 1}$.
\end{cor}

\section{CONCLUSIONS}
It is hoped that the new mathematical ideas presented in this paper would serve as a ``spark plug'' for future number theorists who would consider doing ``serious research'' on OPNs and would pave the way for the eventual resolution of the OPN Problem.

\section*{ACKNOWLEDGMENT}
The author would like to thank his thesis adviser, \emph{Dr. Severino V. Gervacio}, whose persistent prodding and encouragement throughout the author's thesis writing stage made the results presented in this paper possible.  The author's profuse thanks also go to his thesis panelists, \emph{Dr. Leonor A. Ruivivar}, \emph{Dr. Fidel R. Nemenzo} and \emph{Mrs. Sonia Y. Tan}, without whose generosity and keen attention to details, the source master's thesis of the results presented in this paper would not have garnered the Outstanding Graduate Thesis Award.

\end{document}